# Multiple Testing and Error Control in Gaussian Graphical Model Selection

**Mathias Drton and Michael D. Perlman**

*Abstract.* Graphical models provide a framework for exploration of multivariate dependence patterns. The connection between graph and statistical model is made by identifying the vertices of the graph with the observed variables and translating the pattern of edges in the graph into a pattern of conditional independences that is imposed on the variables' joint distribution. Focusing on Gaussian models, we review classical graphical models. For these models the defining conditional independences are equivalent to vanishing of certain (partial) correlation coefficients associated with individual edges that are absent from the graph. Hence, Gaussian graphical model selection can be performed by multiple testing of hypotheses about vanishing (partial) correlation coefficients. We show and exemplify how this approach allows one to perform model selection while controlling error rates for incorrect edge inclusion.

## 1. INTRODUCTION

Many models from multivariate statistics are specified by combining hypotheses of (conditional) independence with particular distributional assumptions. In order to represent such models in a way that is easy to visualize and communicate, it is natural to draw a graph with one vertex for each variable and an edge between any two variables that exhibit a desired type of dependence. In graphical modeling (Cox and Wermuth (1996), Edwards (2000), Lauritzen (1996), Whittaker (1990), Studený (2005)), a rigorous version of this idea is used to associate a statistical model with a graph. Via so-called Markov properties, the pattern of edges in the graph is translated into conditional independence statements, which are then imposed on the joint distribution of the variables that are identified with the graph's vertices. In this process, different graphs with different types of edges have been equipped with different Markov properties. For example, graphs with undirected edges have been given an interpretation that requires two variables that are not joined by an edge to be conditionally independent given all other variables. Markov chains, Markov random fields and certain types of hierarchical log-linear models are examples of models that can be represented in this way. Graphs with directed edges have been used to encode dependence structures that arise from cause-effect relationships among variables (Lauritzen (2001), Pearl (2000), Spirtes, Glymour and Scheines (2000)), and the associated directed graphical models are also known as Bayesian networks. Other types of graphs, sometimes featuring different types of edges simultaneously, have been used to represent other dependence structures. We note that in graphs featuring directed edges, directed cycles, which at an intuitive level correspond to feedback loops, are typically forbidden.

*Mathias Drton is Assistant Professor, Department of Statistics, University of Chicago, Chicago, Illinois 60637, USA e-mail: drton@galton.uchicago.edu. Michael D. Perlman is Professor, Department of Statistics, University of Washington, Seattle 98195, Washington, USA e-mail: michael@ms.washington.edu.*

Much of the success of graphical models in applications, such as the classic application in probabilistic expert systems (Cowell et al. (1999)), is due to favorable computational properties. The independences imposed on the distributions in a graphical model typically induce factorizations of joint densities into smaller, more tractable pieces. The graphical representation of the model then helps to organize computations with these pieces in order to solve statistical inference problems efficiently (Jordan (2004)). In fact, the models are sometimes defined as families of distributions with densities factoring according to a given graph (see, e.g., Jensen (2001)). Conditional independences are then viewed as consequences of such density factorization.

Recently, graphical models have been applied frequently to the analysis of biological data (see, e.g., Beerenwinkel and Drton (2007), Jojic et al. (2004), Lauritzen and Sheehan (2003), McAuliffe, Pachter and Jordan (2004)). In particular, the abundance of gene expression data from microarray experiments has stimulated work on exploratory data analysis focusing on model selection. In the graphical context, this amounts to selection of the underlying graph which may reveal aspects of the network regulating the expression of the genes under study; see Butte et al. (2000), Castelo and Roverato (2006), Dobra et al. (2004), Magwene and Kim (2004), Li and Gui (2006), de la Fuente et al. (2004), Matsuno et al. (2006), Wille et al. (2004), Schäfer and Strimmer (2005) and the review by Friedman (2004).

Three approaches to graphical model selection are commonly taken. The constraint-based approach, which has a long history, is the simplest and employs statistical tests of the model-defining conditional independence hypotheses (Wermuth (1976), Badsberg (1992); Edwards and Havránek (1985), 1987; Kreiner (1987), Smith (1992), Spirtes, Glymour and Scheines (2000), Drton and Perlman (2004)). A second method is a score-based search in which models are selected by searching through the space of underlying graphs and maximizing a goodness-of-fit score such as the Bayesian Information Criterion (BIC) (Schwarz (1978)). The search is often done greedily by defining a neighborhood structure for graphs and terminating with a graph for which no neighboring graph achieves a higher score. While moving in the space of undirected graphs is straightforward by single-edge additions and deletions, this is less simple for graphs with directed edges due to the acyclicity conditions that are usually imposed (see, e.g., Chickering (2002)). Finally, the methodologically most demanding approach to model selection is the Bayesian approach. It requires specification of appropriate prior distributions (Dawid and Lauritzen (1993), Roverato (2002), Roverato and Consonni (2004), Atay-Kayis and Massam (2005)) and computation of/sampling from the resulting posterior distribution on the space of models. Bayesian model determination has been studied for undirected and directed graphs (Cooper and Herskovits (1992), Madigan and Raftery (1994), Heckerman, Geiger and Chickering (1995),Giudici and Green (1999), Consonni and Leucari (2001), Dellaportas, Giudici and Roberts (2003)).

In this paper we consider Gaussian graphical models, which are used in particular for analysis of the continuous gene expression measurements. In Section 2 we review Gaussian graphical models based on undirected, bidirected and acyclic directed graphs. The latter graphs are also known as acyclic digraphs, directed acyclic graphs or, in short, DAGs. We also comment on other graphs that have been used to represent statistical models. Many of these induce Gaussian models that are fully specified by a pairwise Markov property that associates one conditional independence statement with each pair of vertices that are nonadjacent, that is, not joined by an edge. For some graphs, such as undirected and bidirected ones, the conditioning set in such a pairwise conditional independence statement does not depend on the structure of the graph, which is important in our subsequent approach to the model selection problem in which the graph is unknown. For other types of graphs, the same may hold only if a priori information is available that allows one to restrict attention to a restricted subset of graphs. For example, for DAGs, this a priori information may take the form of a total order among the variables, which determines the orientation of any directed edge that is deemed to be present in the graph.

When the absence of edges corresponds to pairwise conditional independence statements, model selection can be performed by testing each conditional independence statement individually. By translating the pattern of rejected hypotheses into a graph, one obtains a constraint-based model selection method in which error rates for incorrect edge inclusion can be controlled when the multiple testing problem is addressed appropriately. Traditionally, such multiple testing approaches are deemed to be not very powerful (Smith (1992)), and instead, classical constraint-based methods use sequential tests in schemes such as, for example, stepwise forward/backward selection. Such sequential procedures may possess greater power for determining the true graph, but the link between the significance level of individual hypothesis tests and overall error properties of the resulting model selection procedure is generally not clear. However, in light of recent advances in multiple testing, it seems worthwhile to revisit the multiple testing approach to graphical model selection. The recent progress not only provides more powerful multiple testing procedures but also allows one to control different types of error rates such as (generalized) family-wise error rate, tail probability of the proportion of false positives and false discovery rate (Sections 3.3 and 3.4). We illustrate the methodology in examples of exploratory data analysis (Section 4), in which the multiple testing approach allows us to identify the most important features of the observed correlation structure. Before concluding in Section 6, we show how prior knowledge about the absence or presence of certain edges can be exploited in order to test fewer and possibly simpler hypotheses in the model selection procedure (Section 5).

## 2. GAUSSIAN GRAPHICAL MODELS

Let $Y = (Y_1, \ldots, Y_p)^t \in \mathbb{R}^p$ be a random vector distributed according to the multivariate normal distribution $\mathcal{N}_p(\mu, \Sigma)$. It is assumed throughout that the covariance matrix $\Sigma$ is nonsingular. Let $G = (V, E)$ be a graph with vertex set $V = \{1, \ldots, p\}$ and edge set $E$. The connection between graph and statistical model is made by identifying the vertices $V$ of the graph $G$ with the variables $Y_1, \ldots, Y_p$. Then the edge set $E$ induces conditional independences via so-called Markov properties. In order to be able to represent different types of dependence patterns, different types of graphs have been equipped with different Markov properties.

### 2.1 Undirected Graphical Models

Let $G = (V, E)$ be an undirected graph, that is, all edges in the graph are undirected edges $i - j$. The *pairwise undirected Markov property* of $G$ associates the conditional independence

$$(2.1) \qquad Y_i \perp\!\!\!\perp Y_j \mid Y_{V \setminus \{i,j\}}$$

with all pairs $(i, j)$, $1 \le i < j \le p$, for which the edge $i - j$ is absent from $G$. For example, the pairwise Markov property for the undirected graph in Figure 1(a) specifies that $Y_1 \perp\!\!\!\perp Y_3 \mid (Y_2, Y_4)$, $Y_1 \perp\!\!\!\perp Y_4 \mid (Y_2, Y_3)$ and $Y_2 \perp\!\!\!\perp Y_3 \mid (Y_1, Y_4)$. The Gaussian graphical model $N(G)$ associated with the undirected graph $G$ is defined as the family of all $p$-variate normal distributions $\mathcal{N}_p(\mu, \Sigma)$ that obey the conditional independence restrictions (2.1) obtained from the pairwise Markov property. Since the random vector $Y$ is distributed according to the multivariate normal distribution $\mathcal{N}_p(\mu, \Sigma)$, we have the equivalence

$$(2.2) \qquad Y_i \perp\!\!\!\perp Y_j \mid Y_{V \setminus \{i,j\}} \quad \Longleftrightarrow \quad \rho_{ij \cdot V \setminus \{i,j\}} = 0,$$

where $\rho_{ij \cdot V \setminus \{i,j\}}$ denotes the $ij$th partial correlation, that is, the correlation between $Y_i$ and $Y_j$ in their conditional distribution given $Y_{V \setminus \{i,j\}}$. This partial correlation can be expressed in terms of the elements of the *concentration* $\equiv$ *precision* matrix $\Sigma^{-1} = \{\sigma^{ij}\}$,

$$(2.3) \qquad \rho_{ij \cdot V \setminus \{i,j\}} = \frac{-\sigma^{ij}}{\sqrt{\sigma^{ii} \sigma^{jj}}};$$

compare Lauritzen (1996, page 130).

The model $N(G)$ has also been called a *covariance selection model* (Dempster (1972)) and a *concentration graph model* (Cox and Wermuth (1996)). The latter name reflects the fact that $N(G)$ can easily be parametrized using the concentration matrix $\Sigma^{-1}$.

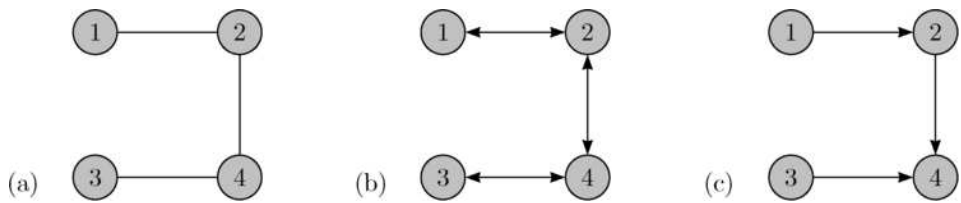

FIG. 1. (a) *An undirected graph*, (b) *a bidirected graph and* (c) *an acyclic directed graph (DAG)*.

Gaussian undirected graphical models are regular exponential families with well-developed statistical methodology (Edwards (2000), Lauritzen (1996), Whittaker (1990)).

Let $Y$ be a random vector whose joint distribution $\mathcal{N}_p(\mu, \Sigma)$ is in the model $N(G)$. By definition the pairwise conditional independences (2.1) hold among the components of $Y$. However, these generally imply other independence relations. For example, if $G$ is the graph from Figure 1(a), then the implied independence relations include

$$(2.4) \qquad Y_1 \perp\!\!\!\perp Y_3 \mid Y_4$$

and

$$(2.5) \qquad Y_1 \perp\!\!\!\perp (Y_3, Y_4) \mid Y_2.$$

One of the benefits of representing the statistical model graphically is that all such independence consequences can be read off the graph using a criterion known as the global Markov property. This criterion is based on paths in the graph, where a path is defined as a sequence of distinct vertices such that any two consecutive vertices in that sequence are joined by an edge. For disjoint subsets $A$, $B$ and $C$ of the vertex set $V$, the *global undirected Markov property* of $G$ states that

$$(2.6) \qquad Y_A \perp\!\!\!\perp Y_B \mid Y_C$$

if there does not exist a path in the graph that leads from a vertex in $A$ to a vertex in $B$ and has *no* nonendpoint vertex in $C$. In other words, the set of vertices $C$ separates the vertices in $A$ from those in $B$. Note that $C$ may be the empty set, in which case conditional independence given $Y_\varnothing$ is understood to be marginal independence of $Y_A$ and $Y_B$. In the graph from Figure 1(a), there is no path from vertex 1 to vertex 3 (or 4) that does not go through vertex 2, which yields (2.5); (2.4) is obtained similarly.

The conditional independence statements (2.1) are saturated in the sense that they involve all the variables at hand. At the other end of the spectrum of pairwise independence statements is marginal independence, the graphical representation of which we discuss next.

## 2.2 Bidirected Graphical Models

Let $G = (V, E)$ be a bidirected graph with edges drawn as $i \longleftrightarrow j$. The *pairwise bidirected Markov property* of $G$ associates the marginal independence

$$(2.7) \qquad Y_i \perp\!\!\!\perp Y_j$$

with all pairs $(i, j)$, $1 \le i < j \le p$, for which the edge $i \longleftrightarrow j$ is absent from $G$. The graph in Figure 1(b), for example, leads to $Y_1 \perp\!\!\!\perp Y_3$, $Y_1 \perp\!\!\!\perp Y_4$ and $Y_2 \perp\!\!\!\perp Y_3$. The Gaussian graphical model $N(G)$ associated with the bidirected graph $G$ is defined as the family of all $p$-variate normal distributions $\mathcal{N}_p(\mu, \Sigma)$ that satisfy the marginal independence restrictions (2.7). Obviously, under multivariate normality,

$$(2.8) \qquad Y_i \perp\!\!\!\perp Y_j \quad \Longleftrightarrow \quad \rho_{ij} = 0,$$

where

$$(2.9) \qquad \rho_{ij} = \frac{\sigma_{ij}}{\sqrt{\sigma_{ii}\sigma_{jj}}}$$

denotes the $ij$th correlation, that is, the correlation between $Y_i$ and $Y_j$.

The model $N(G)$ has also been called a *covariance graph model* (Cox and Wermuth (1996)). We note that Cox and Wermuth (1993, 1996) and some other authors have used dashed instead of bidirected edges. Gaussian bidirected graphical models are curved exponential families, and the development of theory and methodology for these models is still in progress (Chaudhuri, Drton and Richardson (2007), Kauermann (1996), Wermuth, Cox and Marchetti (2006)).

The duality between saturated pairwise conditional independence and marginal independence leads to a nice duality between the global Markov properties for undirected and bidirected Gaussian graphical models. Let $A$, $B$ and $C$ be disjoint subsets of the vertex set $V$. The *global bidirected Markov property* for the graph $G$ states that

$$(2.10) \qquad Y_A \perp\!\!\!\perp Y_B \mid Y_C$$

if there does not exist a path in the graph that leads from a vertex in $A$ to a vertex in $B$ and has *every* nonendpoint vertex on the path in $C$. In the graph from Figure 1(b), it holds that $Y_1 \perp\!\!\!\perp Y_3 \mid Y_2$ because the unique path from vertex 1 to vertex 3 contains the vertex 4 that is not in the conditioning set $\{2\}$. For background regarding the Markov properties of bidirected graphs see Pearl and Wermuth (1994), Kauermann (1996), Banerjee and Richardson (2003) and Richardson (2003).

Both undirected and bidirected graphs have edges without directionality, and any two vertices are either joined by an edge or not. Consequently, in each case there exists a unique complete graph, that is, a graph in which all vertices are joined by an edge.

Moreover, the conditioning sets in the pairwise conditional independences (2.2) and (2.8) do not depend on the structure of the graph. For the acyclic directed graphs introduced next this is no longer true. Since a directed edge between two vertices $i$ and $j$ may be either $i \longrightarrow j$ or $i \longleftarrow j$, there no longer exists a unique complete graph, and the conditioning set in pairwise independence statements associated with missing edges will depend on certain higher-order aspects of the graph.

### 2.3 Directed Graphical Models

A graph with directed edges $i \longrightarrow j$ is called *acyclic* if it contains no directed cycles. A directed cycle is a path of the form $i \longrightarrow \cdots \longrightarrow i$. In the literature, the term directed acyclic graph is often used to refer to these graphs. While this is somewhat imprecise terminology, it yields the popular acronym "DAG," which we also use here.

Let $G = (V, E)$ be an acyclic directed graph/DAG. The directed edges in $E$ define a partial ordering $\preccurlyeq$ of the vertices $V = \{1, \ldots, p\}$ in which $i \preccurlyeq j$ if $i = j$ or there is a directed path $i \longrightarrow \cdots \longrightarrow j$ from $i$ to $j$ in $G$. Not all pairs of vertices must be comparable with respect to this partial ordering. For example, in the graph in Figure 1(c), vertices 2 and 3 are incomparable. However, the partial order can always be extended (possibly nonuniquely) to a total order $\preccurlyeq$ under which $i \leq j$ whenever $i \preccurlyeq j$. (This may require renumbering the vertices.) Such a numbering is called a *well-numbering* or *topological ordering* of the vertex set $V$. In the example, the vertex set is well-numbered, but exchanging vertex numbers 2 and 3 also yields a well-numbering. In the sequel, we assume that the vertex set $V$ is well-numbered.

The *well-numbered pairwise directed Markov property* of $G$ associates the conditional independence

$$(2.11) \qquad Y_i \perp\!\!\!\perp Y_j \mid Y_{\{1,\ldots,j\}\setminus\{i,j\}}$$

with all pairs $(i, j)$, $1 \leq i < j \leq p$, for which the edge $i \longrightarrow j$ is absent from $G$; compare, for example, (7.2) in Edwards (2000) and (2.5) in Drton and Perlman (2007). Note that the well-numbering and the assumed order $i < j$ preclude the existence of the edge $i \longleftarrow j$. In the example of Figure 1(c), (2.11) specifies that $Y_1 \perp\!\!\!\perp Y_3 \mid Y_2$, $Y_2 \perp\!\!\!\perp Y_3 \mid Y_1$ and $Y_1 \perp\!\!\!\perp Y_4 \mid (Y_2, Y_3)$. The Gaussian graphical model $N(G)$ associated with the DAG $G$ is defined as the family of all $p$-variate normal distributions $\mathcal{N}_p(\mu, \Sigma)$ that obey the restrictions (2.11).

The conditional independences (2.11) defining the model $N(G)$ associated with a DAG $G$ may at first sight appear to have a less clear interpretation than the ones associated with undirected or bidirected graphs. However, perhaps even the contrary is true, as the model $N(G)$ based on the DAG $G$ exhibits a dependence structure that can be expected if the (directed) edges in $G$ represent cause-effect relationships. Thinking in this causal fashion, if there is no edge between vertices $i$ and $j$, then $Y_i$ is not an immediate cause of $Y_j$. Hence, if we condition on the variables $Y_{\{1,\ldots,j\}\setminus\{i,j\}}$, which include all immediate (or direct) causes of $Y_j$, then $Y_i$ should have no effect on $Y_j$ as stated in (2.11).

Since we assume that $Y$ follows a multivariate normal distribution $\mathcal{N}_p(\mu, \Sigma)$, it holds that

$$(2.12) \quad \begin{aligned} &Y_i \perp\!\!\!\perp Y_j \mid Y_{\{1,\ldots,j\}\setminus\{i,j\}} \\ &\iff \quad \rho_{ij\cdot\{1,\ldots,j\}\setminus\{i,j\}} = 0, \end{aligned}$$

where for $C \subseteq V \setminus \{i, j\}$ we define $\rho_{ij\cdot C}$ to be the partial correlation of $Y_i$ and $Y_j$ given $Y_C$. Clearly $\rho_{ij\cdot C}$ is a function of the $(C \cup \{i, j\}) \times (C \cup \{i, j\})$ submatrix of $\Sigma$; compare (2.3).

The model $N(G)$ can be shown to correspond to a system of linear regressions. For each variable $Y_i$ there is one linear regression in which $Y_i$ is the response variable and variables $Y_j$ with $j \longrightarrow i$ in $G$ are the covariates (Wermuth (1980), Andersson and Perlman (1998)). It can be shown that the set of covariance matrices giving rise to distributions in $N(G)$ can be parametrized in terms of regression coefficients and residual variances, which constitute the factors of a Choleski decomposition of $\Sigma^{-1}$.

The definition of $N(G)$ does *not* depend on the choice of the underlying well-numbering, which is not unique in general. This follows because a multivariate normal distribution exhibits the conditional independences (2.11) stated by the well-numbered pairwise directed Markov property if and only if it obeys the conditional independences stated by the more exhaustive global directed Markov property, which does not depend on the choice of the well-numbering; see, for example, Cowell et al. (1999, Theorem 5.14), and Drton and Perlman (2007, Appendix A, Theorem 3). The global directed Markov property for a DAG can be formulated in two equivalent ways (Lauritzen (1996)). One way uses a connection to undirected graphs and their Markov properties via what is known as moralization, which involves path separation in certain augmented subgraphs of the DAG. The other, which we detail next, is based on the so-called d-separation criterion which can be applied in the original DAG, but involves an extended definition of separation in the DAG.

Consider a path consisting of the vertices $i_0, i_1, \ldots, i_m$. A nonendpoint vertex $i_k$, $1 \le k \le m-1$, on a path is said to be a *collider* if the preceding and succeeding edges both have an arrowhead at $i_k$, that is, the configuration $i_{k-1} \longrightarrow i_k \longleftarrow i_{k+1}$ occurs in the path. A nonendpoint vertex which is not a collider is said to be a *noncollider*. Two vertices $i$ and $j$ are *d-separated* given a (possibly empty) set $C$ ($i, j \notin C$) if every path between $i$ and $j$ is *blocked relative to $C$*, that is, there is either (i) a noncollider in $C$, or (ii) a collider $c \notin C$ such that there is no vertex $\bar{c} \in C$ with $c \longrightarrow \cdots \longrightarrow \bar{c}$ in $G$. If $A$, $B$ and $C$ are disjoint subsets of the vertex set $V$, then $A$ and $B$ are said to be d-separated given $C$ if every pair of vertices $(i, j)$ with $i \in A$ and $j \in B$ is d-separated given $C$. The *global directed Markov property* of the DAG $G$ states that

$$\text{(2.13)} \qquad Y_A \perp\!\!\!\perp Y_B \mid Y_C$$

whenever $A$ and $B$ are d-separated given $C$.

In the graph in Figure 1(c), the global Markov property implies that $(Y_1, Y_2) \perp\!\!\!\perp Y_3$ because the unique path from vertex 2 to vertex 3 contains the vertex 4, which is a collider that trivially satisfies property (ii) above because the conditioning set $C$ is empty. In this example, it also holds that $Y_1 \perp\!\!\!\perp (Y_3, Y_4) \mid Y_2$ because vertex 2 is a noncollider on the unique path from vertex 1 to vertex 3 (and 4, resp.), and the conditioning set contains this noncollider.

A difficulty in working with DAGs is the fact that two different DAGs can induce the same statistical model, in which case the two graphs are called Markov equivalent. A characterization of this equivalence can be found in Andersson, Madigan and Perlman (1997). However, if two different DAGs share a well-numbering, then they must induce different statistical models, which can be derived, for example, from Lemma 3.2 in Andersson, Madigan and Perlman (1997).

### 2.4 Related Models With Graphical Representation

Many other models considered in the literature benefit from a graphical representation. Most closely connected to the models discussed above are chain graph models. A chain graph is a hybrid graph featuring both undirected and directed edges. The name "chain graph" reflects the fact that the vertex set of these graphs can be partitioned into ordered blocks such that edges between vertices in the same block are undirected, and edges between vertices in different blocks are directed, pointing from the lower-ordered block to the higher-ordered block.

Two alternative Markov properties for chain graphs have been thoroughly studied in the literature: the LWF Markov property of Lauritzen and Wermuth (1989) and Frydenberg (1990), and the more recent AMP Markov property of Andersson, Madigan and Perlman (2001). Results on Markov equivalence and generalizations of the d-separation criterion can be found in Studený and Bouckaert (1998), Studený and Roverato (2006), Andersson and Perlman (2006) and Levitz, Perlman and Madigan (2001). Statistical inference for Gaussian chain graph models is discussed, for example, in Lauritzen (1996) and Drton and Eichler (2006). We note that chain graphs with bidirected instead of undirected edges can also be considered; see, for example, Wermuth and Cox (2004) where dashed edges are used in place of bidirected ones.

Another class of graphical models are based on the ancestral graphs of Richardson and Spirtes (2002). These graphs may include undirected, directed and bidirected edges (although the undirected edges do not occur in an essential way) and permit one to represent all independence structures that may arise from a DAG model under conditioning and marginalization. Ancestral graphs are also related to path diagram/structural equation models that are popular in econometrics and the social sciences; see, for example, Koster (1999). Finally, graphs have also been used to represent time series and stochastic process models (Dahlhaus (2000), Dahlhaus and Eichler (2003), Eichler (2007), Fried and Didelez (2003), Didelez (2007)).

## 3. MODEL SELECTION BY MULTIPLE TESTING

Let $Y^{(1)}, \ldots, Y^{(n)}$ be a sample from a multivariate normal distribution $\mathcal{N}_p(\mu, \Sigma)$ in a Gaussian graphical model $N(G)$, where $G = (V, E)$ is an unknown undirected, bidirected or acyclic directed graph. The sample information can be summarized by the sufficient statistics, which are the sample mean vector

$$\text{(3.1)} \qquad \bar{Y} = \frac{1}{n} \sum_{m=1}^{n} Y^{(m)} \in \mathbb{R}^V$$

and the sample covariance matrix

$$(3.2) \qquad S = \frac{1}{n-1} \sum_{m=1}^{n} (Y^{(m)} - \bar{Y})(Y^{(m)} - \bar{Y})^t$$
$$\in \mathbb{R}^{V \times V}.$$

The problem we consider here is the recovery of the unknown graph underlying the assumed Gaussian graphical model. This is a problem of model selection. We note that in this paper we consider the case where the sample size is moderate to large compared to the number of variables. More precisely, we assume that $n \geq p + 1$ in order to guarantee (almost sure) positive definiteness of the sample covariance matrix $S$. For work on problems in which the sample size is small compared to the number of variables, see, for example, Jones et al. (2005) or Meinshausen and Bühlmann (2006), where sparsity restrictions are imposed on the unknown graph.

## 3.1 Model Selection and Hypotheses of Vanishing Partial Correlations

As presented in Section 2, Gaussian graphical models can be defined by pairwise conditional independence hypotheses or equivalently by vanishing of partial correlations. This suggests that we can perform model selection, that is, recover the graph $G$, by considering the $p(p-1)/2$ testing problems

$$(3.3) \quad H_{ij} : \rho_{ij \cdot C(i,j)} = 0 \quad \text{vs.} \quad K_{ij} : \rho_{ij \cdot C(i,j)} \neq 0$$
$$(1 \leq i < j \leq p).$$

For undirected graphs, (2.2) dictates choosing $C(i, j) = V \setminus \{i, j\}$, and for bidirected graphs we choose $C(i, j) = \varnothing$ in accordance with (2.8). For DAGs, (2.12) leads to the choice $C(i, j) = \{1, \ldots, j\} \setminus \{i, j\}$, which, however, is valid only under the assumption that the vertex set $V = \{1, \ldots, p\}$ is well-numbered for the unknown true DAG.

Thus in order to be able to select a DAG via the testing problems (3.3), we must restrict attention to the situation where we have a priori information about a well-numbering of the vertex set of the unknown DAG. We then select a graph from the set of DAGs for which the specified numbering of the variables is a well-numbering. In an application, a priori information about temporal or causal orderings of the variables can yield a known well-numbering. For example, Spirtes, Glymour and Scheines (2000, Example 5.8.1) analyze data on publishing productivity among academics, which involve seven variables that obey a clear temporal order. In many other applications such a total order among the variables may not be available. However, as we mention in Section 6, a partial order is sufficient for selection of a chain graph (recall Section 2.4).

If in the true graph $G$ there is an edge between vertices $i$ and $j$, then hypothesis $H_{ij}$ is false and the alternative $K_{ij}$ is true. Consequently, if we have performed the $p(p-1)/2$ tests of the hypotheses in (3.3), then we can select a graph by drawing an edge between $i$ and $j$ if and only if the hypothesis $H_{ij}$ is rejected. Let $\alpha \in (0, 1)$ be the significance level employed, and let $\pi_{ij}$ be the $p$-value of the test of hypothesis $H_{ij}$ in (3.3). Then the graph $\hat{G}(\alpha)$ that is selected at level $\alpha$ has the adjacency matrix $\hat{A}(\alpha) = (\hat{a}_{ij}(\alpha)) \in \mathbb{R}^{p \times p}$ with entries

$$(3.4) \qquad \hat{a}_{ij}(\alpha) = \begin{cases} 1, & \text{if } \pi_{ij} \leq \alpha, \\ 0, & \text{if } \pi_{ij} > \alpha. \end{cases}$$

In the sequel we will focus on addressing the issue of multiple testing in this approach, which leads to model selection procedures in which overall error rates (with respect to false inclusion of edges) can be controlled.

We remark that testing the hypotheses in (3.3) is also the first step in stepwise model selection procedures (Edwards (2000), Section 6.1; also see Section 3.1). In *backward stepwise selection*, for example, each hypothesis in (3.3) is tested individually at a fixed significance level $\alpha$. The largest of the $p$-values for the hypotheses that are not rejected is determined and the associated edge is removed from the graph. In the next step the remaining edges/hypotheses are tested again in the reduced graph, also at level $\alpha$. The procedure stops if all remaining hypotheses are rejected at level $\alpha$. While retesting in a reduced graph allows one to take advantage of sparsity of the graph, which induces independence features and may allow for more efficient parameter estimation and testing, the "overall error properties [of such stepwise selection procedures] are not related in any clear way to the error levels of the individual tests" (Edwards (2000), page 158).

## 3.2 Sample Partial Correlations

A natural test statistic for testing hypothesis $H_{ij} : \rho_{ij \cdot C(i,j)} = 0$ is the *sample partial correlation* $r_{ij \cdot C(i,j)}$, that is, the partial correlation computed from the sample covariance matrix $S$; recall (2.3). The marginal distribution of $r_{ij \cdot C(i,j)}$ has the same form as the distribution of the ordinary sample correlation $r_{ij}$, but with the parameter $\rho_{ij}$ replaced by $\rho_{ij \cdot C(i,j)}$ and the sample size $n$ reduced to $n_{C(i,j)} = n - |C(i,j)|$ (Anderson (2003), Theorem 4.3.5). The marginal distribution of the sample correlation $r_{ij}$ takes on a simple form if the $i$th and $j$th components of the normal random vector from which it is derived are independent.

PROPOSITION 1. *If the true correlation $\rho_{ij}$ is zero, then $\sqrt{n-2} \cdot r_{ij} / \sqrt{1 - r_{ij}^2}$ has a $t$-distribution with $n-2$ degrees of freedom.*

In the noncentral case, $\rho_{ij} \neq 0$, the exact distribution of $r_{ij}$ can be described using hypergeometric functions, but it is simpler to work with Fisher's variance-stabilizing $z$-transform.

PROPOSITION 2. *Let*

$$z : (-1, 1) \to \mathbb{R}, \quad r \mapsto \frac{1}{2} \ln \left( \frac{1+r}{1-r} \right)$$

*be the $z$-transform. An accurate normal approximation to the distribution of $z_{ij} = z(r_{ij})$ can be obtained from the fact that*

$$\sqrt{n-3} \, (z_{ij} - \zeta_{ij}) \xrightarrow{d} \mathcal{N}(0, 1) \quad \text{as } n \to \infty,$$

*where $\zeta_{ij} = z(\rho_{ij})$. Note that $\zeta_{ij} = 0$ if and only if $\rho_{ij} = 0$.*

Concerning finite-sample properties, little is lost by working with the normal approximation to Fisher's $z$, as this approximation is accurate even for moderate sample size (Anderson (2003), Section 4.2.3). At the same time much convenience is gained because of the variance-stabilizing property of the $z$-transform and the fact that the joint distribution of $z$-transformed sample correlations is easily deduced from that of the untransformed sample correlations; compare Proposition 5 below.

A sample partial correlation $r_{ij \cdot C(i,j)}$ is a smooth function of the sample covariance matrix $S$. The random matrix $(n-1)S$ has a Wishart distribution with $n-1$ degrees of freedom and scale parameter matrix $\Sigma$. The asymptotic normal distribution of both $S$ and $S^{-1}$ can be described using Isserlis matrices (Olkin and Siotani (1976), Roverato and Whittaker (1998)).

PROPOSITION 3. *Let $\text{Iss}(\Sigma)$ be the Isserlis matrix of $\Sigma$, that is, the $p(p+1)/2 \times p(p+1)/2$-matrix with entries*

$$\text{Iss}(\Sigma)_{ij,uv} = \sigma_{iu}\sigma_{jv} + \sigma_{iv}\sigma_{ju},$$
$$1 \leq i \leq j \leq p, 1 \leq u \leq v \leq p.$$

*Then*

$$\sqrt{n}(S - \Sigma) \xrightarrow[d]{n \to \infty} \mathcal{N}_{p(p+1)/2}(0, \text{Iss}(\Sigma)),$$

*and*

$$\sqrt{n}(S^{-1} - \Sigma^{-1}) \xrightarrow[d]{n \to \infty} \mathcal{N}_{p(p+1)/2}(0, \text{Iss}(\Sigma^{-1})).$$

Using the delta method (van der Vaart (1998)), the joint asymptotic normal distribution of the vector of sample partial correlations $r = (r_{ij \cdot C(i,j)} \mid 1 \leq i < j \leq p)$ can be derived. For ordinary correlations, for which $C(i,j) = \varnothing$ for all $1 \leq i < j \leq p$, and saturated partial correlations with $C(i,j) = V \setminus \{i, j\}$ for all $1 \leq i < j \leq p$, the following result is quickly obtained using software for symbolic computation. The statement about ordinary correlations goes back to Aitkin (1969, 1971) and Olkin and Siotani (1976).

PROPOSITION 4. *The vector of ordinary correlations is asymptotically normal,*

$$\sqrt{n}(r - \rho) \xrightarrow[d]{n \to \infty} \mathcal{N}_{p(p-1)/2}(0, \Omega),$$

*with $\rho = (\rho_{ij} \mid 1 \leq i < j \leq p)$ and the asymptotic covariance matrix $\Omega = (\omega_{ij,k\ell})$ given by*

$$\omega_{ij,ij} = [1 - (\rho_{ij})^2]^2,$$
$$\omega_{ij,i\ell} = -\frac{1}{2}\rho_{ij}\rho_{i\ell}[1 - (\rho_{ij})^2 - (\rho_{i\ell})^2 - (\rho_{j\ell})^2]$$
$$\quad + \rho_{j\ell}[1 - (\rho_{ij})^2 - (\rho_{i\ell})^2],$$
$$\omega_{ij,k\ell} = \frac{1}{2}\rho_{ij}\rho_{k\ell}[(\rho_{ik})^2 + (\rho_{i\ell})^2 + (\rho_{jk})^2 + (\rho_{j\ell})^2]$$
$$\quad + \rho_{ik}\rho_{j\ell} + \rho_{i\ell}\rho_{jk} - \rho_{ik}\rho_{jk}\rho_{k\ell}$$
$$\quad - \rho_{ij}\rho_{ik}\rho_{i\ell} - \rho_{ij}\rho_{jk}\rho_{j\ell} - \rho_{i\ell}\rho_{j\ell}\rho_{k\ell}.$$

*The same result holds for the vector of saturated partial correlations $r_{ij \cdot V \setminus \{i,j\}}$ if we replace all $\rho_{ij}$ by $\rho_{ij \cdot V \setminus \{i,j\}}$ in the above formulas, where, for more accurate normal approximation, the sample size $n$ should also be replaced by $n_{V \setminus \{i,j\}} = n - p - 2$.*

For small number of variables $p$, asymptotic covariance matrices for vectors of partial correlations $r_{ij \cdot \{1, \ldots, j\} \setminus \{i,j\}}$, as required for DAG selection, can be computed using software for symbolic computation, but we are not aware of any general formulas in the literature.

By again applying the delta method, the following result is obtained:

PROPOSITION 5. *The asymptotic covariance matrix of the vector of $z$-transformed partial correlations $z_{ij \cdot C(i,j)} = z(r_{ij \cdot C(i,j)})$ is the correlation matrix of the asymptotic covariance matrix $\Omega$ of the vector of untransformed partial correlations $r_{ij \cdot C(i,j)}$.*

With these preliminaries we can now turn to the problem of error control in the multiple testing problem (3.3), which we formulate in terms of $p$-values.

### 3.3 Controlling Family-Wise Error Rate

A multiple testing procedure for problem (3.3) is said to control the family-wise error rate (FWER) at level $\alpha \in (0, 1)$ if for any underlying multivariate normal distribution $\mathcal{N}_p(\mu, \Sigma)$, the probability of rejecting one or more null hypotheses $H_{ij}$ incorrectly is smaller than or equal to $\alpha$. If a multiple testing procedure controls the FWER at level $\alpha$, then its (simultaneous) $p$-values $\{\pi_{ij} \mid 1 \leq i < j \leq p\}$ have the property

$$\text{Prob}_{\mathcal{N}_p(\mu, \Sigma)}(\exists_{ij} : \pi_{ij} \leq \alpha \text{ but } H_{ij} \text{ is true})$$
$$= \text{Prob}_{\mathcal{N}_p(\mu, \Sigma)}(\exists_{ij} : \text{edge } i \text{ --- } j \text{ included}$$
$$\text{when actually absent}) \leq \alpha.$$

Since control is achieved at all multivariate normal distributions, that is, at all patterns of true null hypotheses, this is sometimes referred to as "strong control" (Dudoit, Shaffer and Boldrick (2003)). For the graph selected according to (3.4), this means that with probability at most $\alpha$ the selected graph $\hat{G}(\alpha)$ is not a subgraph of the true graph $G$,

$$(3.5) \qquad \text{Prob}_G(\hat{G}(\alpha) \nsubseteq G) \leq \alpha.$$

(By definition, a subgraph of $G$ contains no edges that are absent in $G$.) The notation $\text{Prob}_G$ in (3.5) denotes any probability calculation under a distribution $\mathcal{N}_p(\mu, \Sigma) \in N(G)$.

In (3.5), the error with respect to incorrect edge inclusion is controlled in finite samples. Sometimes it may, however, only be feasible to achieve asymptotic control of the form

$$\limsup_{n \to \infty} \text{Prob}_{\mathcal{N}_p(\mu, \Sigma)}(\exists_{ij} : \pi_{ij} \leq \alpha \text{ but } H_{ij} \text{ is true}) \leq \alpha,$$

and hence,

$$(3.6) \qquad \limsup_{n \to \infty} \text{Prob}_G(\hat{G}(\alpha) \nsubseteq G) \leq \alpha.$$

In particular, when using $z$-transformed correlations and normal approximations to their distributions, one may only hope to achieve asymptotic control given by (3.6). But the concept of asymptotic control is also central to recently introduced multiple testing procedures that employ consistent estimates of the joint distribution of the test statistics (see Section 3.3.2).

Any model selection procedure that asymptotically controls FWER at fixed significance level $\alpha$, that is, satisfies (3.6), is $(1 - \alpha)$-consistent in the following sense. Let $G_{\text{faithful}} = G_{\text{faithful}}(\Sigma)$ be the graph that has the fewest edges among all graphs $G'$ for which the data-generating distribution $\mathcal{N}_p(\mu, \Sigma)$ is in $N(G')$. The definition of $G_{\text{faithful}}$ is straightforward for undirected and bidirected graphs. For DAGs, recall that we assume to know a well-numbering of the variables a priori, in which case $G_{\text{faithful}}$ is again well defined. The distribution $\mathcal{N}_p(\mu, \Sigma)$ is pairwise faithful to $G_{\text{faithful}}$ in the sense that $\rho_{ij \cdot C(i,j)} = \rho_{ij \cdot C(i,j)}(\Sigma) = 0$ if and only if the edge between vertices $i$ and $j$ is absent from $G_{\text{faithful}}$. Then it can be shown that

$$(3.7) \qquad \liminf_{n \to \infty} \text{Prob}_{\mathcal{N}_p(\mu, \Sigma)}(\hat{G}(\alpha) = G_{\text{faithful}}) \geq 1 - \alpha;$$

compare Drton and Perlman [2004, (2.18)]. This means that the selection procedure identifies $G_{\text{faithful}}$ with asymptotic probability at least $(1 - \alpha)$. Moreover, if the sample size $n$ can be chosen large enough, then the asymptotic probability $\text{Prob}_{\mathcal{N}_p(\mu, \Sigma)}(\hat{G}(\alpha) \neq G_{\text{faithful}})$ can be made arbitrarily small by choosing $\alpha$ arbitrarily small [Drton and Perlman (2004), (2.19)–(2.21)]. In this sense the procedure is fully consistent. However, the choice of $n$ depends on $\alpha$ as well as on $\Sigma$.

The above notion of faithfulness is defined with respect to the pairwise Markov property. In other contexts (e.g., Spirtes, Glymour and Scheines (2000), Wille and Bühlmann (2006), Becker, Geiger and Meek (2000)), the stronger condition of faithfulness with respect to the global Markov property is considered. A distribution is globally faithful to a graph $G$ if it exhibits a conditional independence $Y_A \perp\!\!\!\perp Y_B \mid Y_C$ if and only if the sets $A$, $B$ and $C$ fulfill the graphical separation property used in the definition of the global Markov property. However, while every multivariate normal distribution is pairwise faithful to some graph in the considered class of graphs, it is easy to see that there exist normal distributions that are not globally faithful to any graph in the class. For example, consider the centered trivariate normal distribution with covariance matrix

$$\Sigma = \begin{pmatrix} 2 & 1 & 0 \\ 1 & 2 & 1 \\ 0 & 1 & 1 \end{pmatrix}$$
$$\implies \Sigma^{-1} = \begin{pmatrix} 1 & -1 & 1 \\ -1 & 2 & -2 \\ 1 & -2 & 3 \end{pmatrix}.$$

Choosing the context of undirected graphs, we see that this distribution is pairwise but not globally faithful to the complete graph, since $\sigma_{13} = 0$ implies $Y_1 \perp\!\!\!\perp Y_3$.

### 3.3.1 Multiple testing procedures using the marginal distributions of the sample correlations.

Classical generally applicable multiple testing procedures are based on the marginal distributions of the test statistics alone. When testing $H_{ij} : \rho_{ij \cdot C(i,j)} = 0$ via the normal approximation to the $z$-transformed sample partial correlation $z_{ij \cdot C(i,j)}$ given in Proposition 2, we obtain the unadjusted $p$-value

$$(3.8) \qquad \pi_{ij} = 2[1 - \Phi(\sqrt{n_{C(i,j)} - 3} \cdot |z_{ij \cdot C(i,j)}|)],$$

where $\Phi$ is the cumulative distribution function of the standard normal distribution $\mathcal{N}(0, 1)$. These $p$-values can now be adjusted to achieve FWER control (3.6) in the selection of the graph $\hat{G}(\alpha)$ defined in (3.4).

The Bonferroni $p$-values

$$(3.9) \quad \pi_{ij}^{\text{Bonf}} = \min\left\{\binom{p}{2}\pi_{ij}, 1\right\}, \quad 1 \le i < j \le p,$$

are the simplest such adjusted $p$-values. An easy, more powerful adjustment is obtained by the step-down method of Holm (1979). In this method the $p$-values $\pi_{ij}$ in (3.8) are ordered as $\pi_{1\uparrow} \le \pi_{2\uparrow} \le \cdots$ and the ordered adjusted $p$-values are obtained as

$$(3.10) \quad \begin{aligned} &\pi_{a\uparrow}^{\text{Bonf,Step}} \\ &= \max_{b=1,\ldots,a}\left[\min\left\{\left(\binom{p}{2} - b + 1\right)\pi_{b\uparrow}, 1\right\}\right], \\ &\hspace{5cm} 1 \le a \le \binom{p}{2}. \end{aligned}$$

The adjusted $p$-value $\pi_{a\uparrow}^{\text{Bonf,Step}}$ is then associated with the hypothesis $H_{ij}$ that gave rise to the $a$th smallest unadjusted $p$-value. Note that if the original $p$-values in (3.8) were computed using the $t$-transform in Proposition 1 instead of Fisher's $z$-transform, then both these $p$-value adjustments would provably achieve finite-sample error control (3.5).

A more powerful procedure than Bonferroni is obtained by applying Šidák's inequality (Šidák (1967)) to the joint asymptotic normal distribution of the vector of transformed sample partial correlation coefficients $z_{ij \cdot C(i,j)}$. The inequality yields the $p$-values

$$(3.11) \quad \pi_{ij}^{\text{Sidak}} = 1 - (1 - \pi_{ij})^{\binom{p}{2}}, \quad 1 \le i < j \le p,$$

where again $\pi_{ij}$ is given by (3.8). The $p$-values $\pi_{ij}^{\text{Sidak}}$, which appear in Drton and Perlman [2004, (2.9)], can in turn be improved in a step-down approach to

$$(3.12) \quad \begin{aligned} \pi_{a\uparrow}^{\text{Sidak,Step}} &= \max_{b=1,\ldots,a}[1 - (1 - \pi_{b\uparrow})^{(\binom{p}{2} - b + 1)}], \\ &\hspace{4cm} 1 \le a \le \binom{p}{2}. \end{aligned}$$

As in the case of the Bonferroni-adjusted $p$-values, the index $a\uparrow$ refers to the $a$th smallest $p$-value and the ordered adjusted $p$-value $\pi_{a\uparrow}^{\text{Sidak,Step}}$ is to be associated with the hypothesis $H_{ij}$ that gave rise to the $a$th smallest unadjusted $p$-value. Both sets of $p$-values, $\pi_{ij}^{\text{Sidak}}$ and $\pi_{ij}^{\text{Sidak,Step}}$, define a graph $\hat{G}(\alpha)$ that satisfies (3.6).

### 3.3.2 Multiple testing procedures using the joint distribution of the sample correlations.

Westfall and Young (1993) describe multiple testing methods that can improve upon the marginal distribution-based procedures from Section 3.3.1 by exploiting possible dependences among the test statistics used to test the individual hypotheses. However, these methods cannot be applied to testing of correlations, as the required condition known as "subset-pivotality" is not satisfied in this context (Westfall and Young (1993), page 43). A way around this condition was found recently by Pollard and van der Laan (2004), Dudoit, van der Laan and Pollard (2004) and van der Laan, Dudoit and Pollard (2004b), who describe how a consistent estimate of the asymptotic joint multivariate normal distribution of the test statistics can indeed be used for a valid $p$-value adjustment. We now detail this approach in our context.

For sample size $n$ tending to infinity, our vector of test statistics $z = (z_{ij \cdot C(i,j)} \mid 1 \le i < j \le p)$ has a multivariate normal limiting distribution which can be derived from that of $r = (r_{ij \cdot C(i,j)} \mid 1 \le i < j \le p)$ using Proposition 5. We obtain the normal approximations

$$(3.13) \quad \begin{aligned} r &\overset{.}{\sim} \mathcal{N}_{p(p-1)/2}(\rho, N^{-1}\Omega N^{-t}) \quad \Longrightarrow \\ z &\overset{.}{\sim} \mathcal{N}_{p(p-1)/2}(\zeta, N^{-1}\,\text{Corr}(\Omega)N^{-t}), \end{aligned}$$

where $\rho = (\rho_{ij \cdot C(i,j)} \mid 1 \le i < j \le p)$, $\zeta$ is the component-wise $z$-transform of $\rho$, and $N$ is the diagonal matrix with diagonal entries $\sqrt{n_{C(i,j)} - 3}$. Recall that for undirected and bidirected graphs with $C(i,j) = V \setminus \{i,j\}$ and $C(i,j) = \varnothing$, respectively, Proposition 3 yields the asymptotic covariance matrix $\Omega$, whereas for DAGs with $C(i,j) = \{1,\ldots,j\} \setminus$

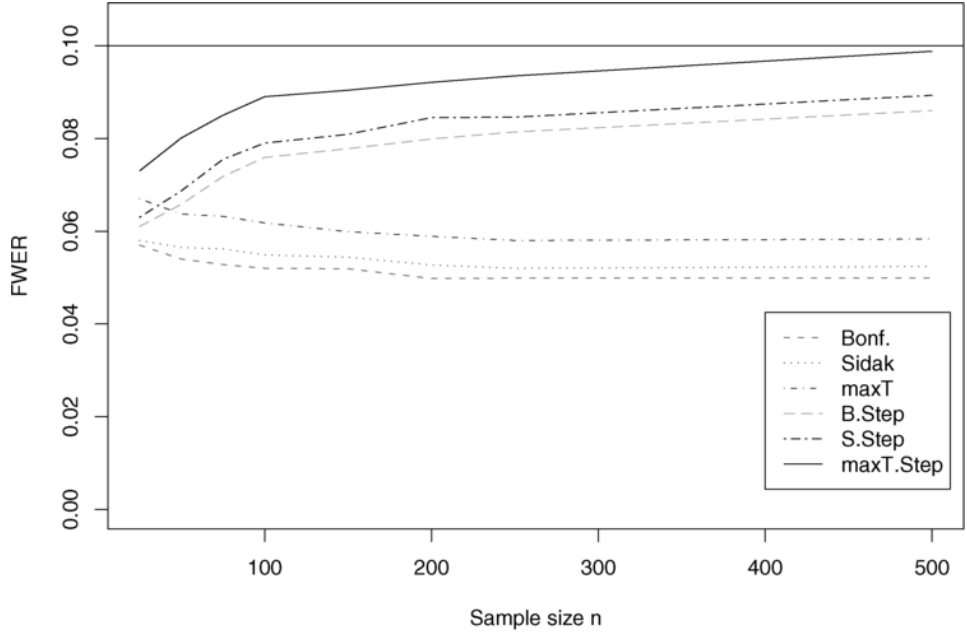

F<small>IG</small>. 2. *Simulated family-wise error rates for false edge inclusion in undirected graph model selection. The multiple testing procedures were set up to control the family-wise error rate at level $\alpha = 0.1$ and 10,000 samples were drawn from a multivariate normal distribution for $p = 7$ variables; sample size varied from $n = 25$ to $n = 500$. In the normal distribution, nine partial correlations were nonzero with values in $[0.2, 0.55]$.*

$\{i, j\}$ no general formula for the asymptotic covariance matrix seems to be available.

The asymptotic covariance matrix $\Omega$ does involve unknown quantities derived from the covariance matrix $\Sigma$ of our observed random vectors. Plugging in the corresponding expression formed from the sample covariance matrix $S$, we obtain a consistent estimator $\hat{\Omega}$. We can then determine the so-called max-$T$ adjusted $p$-values

$$
\pi_{ij}^{\max} = \text{Prob}_{\mathcal{N}(0, N^{-1}\text{Corr}(\hat{\Omega})N^{-t})}
$$
$$
(3.14) \qquad \left( \max_{1 \le u < v \le p} |Z_{uv}| \ge |z_{ij \cdot C(i,j)}| \right),
$$
$$
1 \le i < j \le p.
$$

These probabilities can be computed by Monte Carlo simulation drawing the vector $Z = (Z_{uv})$ from $\mathcal{N}(0, N^{-1}\text{Corr}(\hat{\Omega})N^{-t})$. A step-down max-$T$ procedure is also available. It is based on ordering the $z_{ij \cdot C(i,j)}$ as $|z_{1\downarrow}| \ge |z_{2\downarrow}| \ge \cdots$ and yields the adjusted $p$-values

$$
\pi_{a\uparrow}^{\max,\text{Step}} = \max_{b=1,\ldots,a} \text{Prob}_{\mathcal{N}(0, N^{-1}\text{Corr}(\hat{\Omega})N^{-t})}
$$
$$
(3.15) \qquad \left( \max_{b=a,\ldots,\binom{p}{2}} |Z_{b\downarrow}| \ge |z_{a\downarrow}| \right),
$$
$$
1 \le a \le \binom{p}{2}.
$$

[Note that since the asymptotic marginal distributions of the $z$-transformed partial correlations $z(r_{ij \cdot C(i,j)})$ are identical, all being standard normal, the min-$P$ adjustment is identical to the max-$T$ adjustment; see Dudoit, Shaffer and Boldrick (2003).]

Both the single-step $p$-values $\pi_{ij}^{\max}$ and the less conservative step-down $p$-values $\pi_{ij}^{\max,\text{Step}}$ define a graph $\hat{G}(\alpha)$ for which the condition (3.6) for asymptotic control of the FWER holds (Pollard and van der Laan (2004); Dudoit, van der Laan and 2004; van der Laan, Dudoit and Pollard 2004b). In fact, it follows from the general results in van der Laan, Dudoit and Pollard (2004b) that the step-down $p$-values $\pi_{ij}^{\max,\text{Step}}$ yield a graph $\hat{G}(\alpha)$ satisfying

$$
(3.16) \qquad \lim_{n \to \infty} \text{Prob}_G(\hat{G}(\alpha) \nsubseteq G) = \alpha.
$$

In other words, the asymptotic error control is not conservative but exact.

The simulations summarized in Figure 2 show that the step-down max-$T$ adjustment method based on (3.15) does indeed provide the most exact error control for false edge inclusion. However, the step-down procedures based on marginal distributions may still be useful if, for large number of variables $p$, the Monte Carlo computation needed to compute the max-$T$ adjusted $p$-values becomes too time-consuming.

### 3.4 Alternative Error Rates

In multiple testing-based graphical model selection, controlling the FWER amounts to controlling the probability that a single one of the edges included in the selected graph is incorrect. This is clearly a very stringent requirement. Alternatively, other less demanding error rates can be controlled, of which we now discuss the three most popular ones; see Dudoit, van der Laan and Pollard (2004), Romano and Wolf (2005) and van der Laan, Dudoit and Pollard (2004a, 2004b) for recent surveys of the relevant literature.

One relaxation consists of controlling the *generalized family-wise error rate* (GFWER), which is defined with respect to a chosen nonnegative integer $k$. This generalization is based on the probability of the event that at most $k$ of the true null hypotheses are incorrectly rejected. In our graphical model selection context, if a multiple testing procedure is set up to control the $k$-GFWER at level $\alpha \in (0,1)$, and its adjusted $p$-values are used to select a graph $\hat{G}(\alpha)$, then it holds that

$$\text{(3.17)} \quad \text{Prob}_G(\hat{G}(\alpha) \text{ contains } k+1 \text{ or more}$$
$$\text{edges that are not present in } G) \leq \alpha.$$

For $k = 0$, (3.17) reduces to (3.5), that is, control of the traditional FWER.

In some contexts, one may be willing to live with a larger number of erroneous edge inclusion decisions if the selected graph is less sparse. This can be achieved by controlling the tail probability of the *proportion of false positives*, also known as *false discovery proportion*. Here a fraction $\lambda \in [0,1)$ is chosen and the probability of the event that more than a proportion $\lambda$ of the rejected hypotheses are incorrectly rejected is to be controlled. In the present context, control of the tail probability of the proportion of false positives (TPPFP) at level $\alpha$ allows us to select a graph $\hat{G}(\alpha)$ with the property that

$$\text{(3.18)} \quad \text{Prob}_G(\text{More than } 100\lambda\% \text{ of the edges}$$
$$\text{in } \hat{G}(\alpha) \text{ are not present in } G) \leq \alpha.$$

For $\lambda = 0$, (3.18) reduces to (3.5).

Of a somewhat different nature is the false discovery rate (FDR), which is defined in terms of the expectation of the proportion of false positives. Controlling the FDR at level $\alpha$ allows us to select a graph $\hat{G}(\alpha)$ such that the proportion of incorrect edges among all the edges of $\hat{G}(\alpha)$ is smaller than $\alpha$ in expectation, that is,

$$\text{(3.19)} \quad \text{E}_G\left[\frac{\#\text{edges incorrectly included in } \hat{G}(\alpha)}{\#\text{edges included in } \hat{G}(\alpha)}\right]$$
$$\leq \alpha.$$

A number of methods for control of GFWER and TPPFP have been described in the literature and can be applied for graphical model selection. Perhaps the simplest methods are the augmentation methods of van der Laan, Dudoit and Pollard (2004a). The idea there is to first determine the hypotheses rejected in FWER control and then reject additional hypotheses. For $k$-GFWER control one simply rejects $k$ additional hypotheses from the most significant not already rejected ones. For $\lambda$-TPPFP control the augmentation method proceeds similarly with the number of additionally rejected hypotheses determined from the parameter $\lambda$ and the number of hypotheses already rejected in FWER control. While simple and asymptotically exact if used in conjunction with the asymptotically exact max-$T$ step-down procedure for FWER control [cf. (3.16)], the augmentation methods may sometimes be outperformed by methods that address the respective generalized error rate directly and not via FWER control; see Romano and Wolf (2005), who survey such methods that can be designed to employ either the joint distribution of the test statistics or only their marginals.

For control of the FDR, the original step-up method of Benjamini and Hochberg (1995) is not generally applicable in our context as it requires the test statistics to exhibit a form of dependence termed "positive regression dependency" (Benjamini and Yekutieli (2001)). A generally valid method is obtained by introducing a log-term as penalty in the step-up method (Benjamini and Yekutieli (2001)). Alternatively, van der Laan, Dudoit and Pollard (2004a) proposed a method for FDR control that is derived from TPPFP control.

## 4. APPLICATION TO GENE EXPRESSION DATA

In this section, we demonstrate multiple testing-based graphical model selection using data from $n = 118$ microarray experiments collected and analyzed by Wille et al. (2004); see also Wille and Bühlmann (2006). The experiments measure gene expression in *Arabidopsis thaliana*, and for our purposes we focus on $p = 13$ genes from the initial part of the MEP pathway, which is one of the two pathways that received special attention in Wille et al. (2004). We select graphical models by multiple testing in order to discover the key features of the correlation structure among the considered gene expression measurements. Revealing such key features is important if the goal of the analysis is to generate scientific hypotheses about the interplay of genes in a gene regulatory network. First, in Section 4.1, we will select different types of graphs with the goal of emphasizing how different graphs capture different types of dependence. Then, in Section 4.2, we give a simple example of the use of alternative error rates.

### 4.1 Selecting Different Graphs

In order to select an undirected graph we apply a multiple testing procedure to the testing problem (3.3) with $C(i,j) = V \setminus \{i,j\}$. Applying the step-down max-$T$ procedure from Section 3.3.2 to control the classical FWER at simultaneous significance level $\alpha = 0.15$, we select the undirected graph depicted in Figure 3, which has 16 edges. Control of the FWER allows us to state that we are 85% confident that all the edges in this graph are also present in the true graph. In this example, the step-down max-$T$ procedure is indeed the most powerful of the procedures described in Sections 3.3.1 and 3.3.2. For example, if $i = $ DXPS1 and $j = $ GPPS, then $\pi_{ij}^{\max, \text{Step}} = 0.066$ but $\pi_{ij}^{\text{Bonf}} = 0.099$.

The undirected graph in Figure 3 shows some of the features of larger graphs shown in Wille et al. (2004), but there are also differences. Note, however, that the graphical modeling approach in Wille et al. (2004) is different from any of the approaches described here in that only conditional independences involving three variables are considered.

Next we select a bidirected graph by testing (3.3) with $C(i,j) = \varnothing$. Using again the step-down max-$T$ procedure with simultaneous significance level $\alpha = 0.15$, we select the bidirected graph in Figure 3. The selected bidirected graph features 30 edges. As for the selection of the undirected graph, the step-down max-$T$ procedure provides the most powerful method for FWER control. Use of any of the other multiple testing procedures from Sections 3.3.1 and 3.3.2 results in the selection of a graph with 28 or 29 edges. The two graphs in Figure 3 have some common edges but many adjacencies are different, which is a reflection of the fact that large correlations are not necessarily associated with large partial correlations and vice versa. The two correlation measures quantify very different types of dependence.

The vertical placement of the vertices in the graphs in Figure 3 reflects a partial order among the considered genes that is based on the genes' role in the metabolic network (Wille et al. (2004), Figure 2). In order to illustrate selection of a DAG, we refine this partial order to a total order in which DXPS1 < DXPS2 < DXPS3 < DXR < $\cdots$ < IPPI1 < GPPS < PPDS1 < PPDS2. We then test the hypotheses (3.3) with $C(i,j) = \{1, \ldots, j\} \setminus \{i,j\}$, where the indices $i$ and $j$ refer to the rank of a gene in the total order. Since the asymptotic covariance matrix of the sample partial correlations used to test these hypotheses is unknown, we use the step-down $p$-values $\pi_{ij}^{\text{Sidak, Step}}$ to control FWER at level $\alpha = 0.15$. (Note, however, that bootstrap-based methods can be used to estimate the unknown joint distribution of the test statistics nonparametrically; compare Dudoit et al. (2004); van der Laan, Dudoit and Pollard, 2004a, 2004b; Romano and Wolf (2005).) The selected DAG, depicted on the left in Figure 4, has 19 edges. Among the three graphs we selected, this DAG best reflects the structure of the metabolic network formed by the considered genes; recall, however, that we used the structure of the metabolic network to form a well-numbering of the variables. We remark that a strict causal interpretation of this DAG (compare Section 2.3) would rest on the assumption that there are no hidden/unobserved causes.

### 4.2 Alternative Error Rates

In order to convey how more liberal error rates allow for the inclusion of additional edges, we consider the selection of DAGs under control of GFWER and TPPFP. Since we have already computed $p$-values for FWER control, the augmentation methods of van der Laan et al. (2004a) can be readily applied.

For control of the $k$-GFWER, we simply determine the $k$ smallest of the FWER $p$-values that are associated with hypotheses not rejected by the FWER controlling procedure. We then reject the $k$ hypotheses corresponding to these $p$-values. Choosing $k = 5$ and keeping the simultaneous significance level $\alpha = 0.15$ used above, we select the DAG with 24 edges shown on the right-hand side in Figure 4. Due to the GFWER control we can state that we are 85% confident that at most five edges in this graph are not present in the true underlying graph.

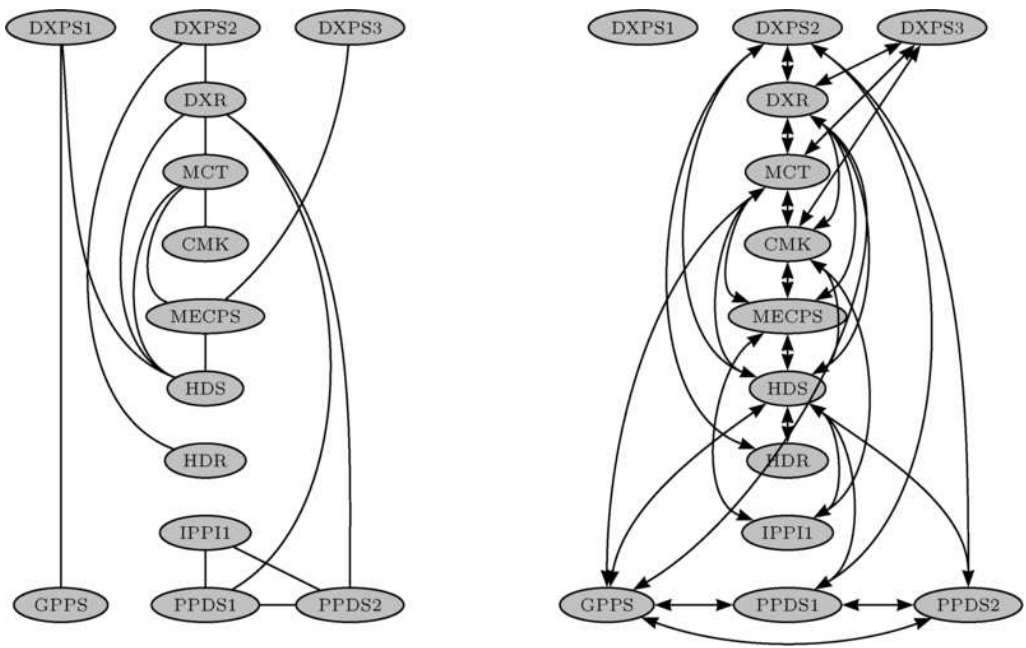

FIG. 3. *An undirected and a bidirected graph selected by controlling FWER at $\alpha = 0.1$ with the step-down max-T procedure.*

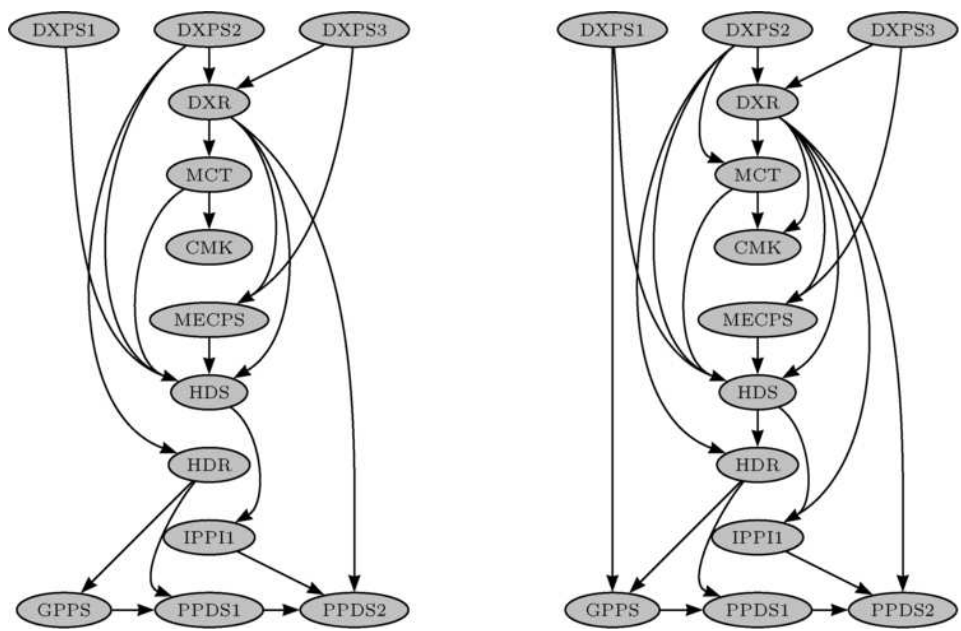

FIG. 4. *Two DAG's selected by the step-down Sidak procedure. The graph to the left is obtained by controlling FWER at level $\alpha = 0.1$, the one to the right by controlling $k$-GFWER with $k = 5$ at $\alpha = 0.1$.*

(We remark that a simple, more direct step-down approach to GFWER control described in Lehmann and Romano (2005), Theorem 2.2, leads to the same DAG.)

TPPFP control by augmentation proceeds again by rejecting additional hypotheses not yet rejected by the procedure for FWER control. Keeping with $\alpha = 0.15$ and choosing the proportion $\lambda = 0.22$, the $\lambda$-TPPFP control by augmentation again yields the graph on the right-hand side in Figure 4. Therefore, we are 85% confident that at most 22% of the edges of this graph are not present in the true underlying DAG.

## 5. INCORPORATING PRIOR INFORMATION ABOUT THE PRESENCE OR ABSENCE OF EDGES

Suppose it is known that in the true graph $G = (V, E)$ certain edges $E_0$ are absent, $E_0 \cap E = \varnothing$, and certain other edges $E_1$ are present, $E_1 \subseteq E$. Model selection then reduces to the problem of determining the absence or presence of the uncertain edges $E_u$, that is, the complement of $E_0 \cup E_1$ in the set of all possible edges. Let $G_{\mathrm{up}} = (V, E_1 \dot\cup E_u)$ denote the *upper graph*, which contains all edges known to be present as well as all uncertain edges. In the context of DAGs, Consonni and Leucari (2001) call the upper graph the "full graph." Similarly, let $G_{\mathrm{low}} = (V, E_1)$ denote the *lower graph*, which includes only the edges that are known to be present. Thus the true graph $G$ satisfies $G_{\mathrm{low}} \subseteq G \subseteq G_{\mathrm{up}}$, where $G_{\mathrm{low}}$ and $G_{\mathrm{up}}$ are known. If all edges are uncertain, then the upper and lower graph are the complete and the empty graph, respectively.

The multiple testing approach presented in Section 3 extends readily to the present case by reducing the $p(p-1)/2$ simultaneous testing problems from (3.3) to the $q = |E_u|$ testing problems corresponding to the uncertain edges only. Since $q \leq p(p-1)/2$, we have fewer testing problems to consider and gain power in simultaneous testing. Furthermore, since $G \subseteq G_{\mathrm{up}}$, the conditional independences holding in $G_{\mathrm{up}}$ also hold in $G$, which may allow for additional power gain because, as we explain next, the hypothesis $H_{ij} : \rho_{ij \cdot C(i,j)} = 0$ may be reformulated equivalently using a smaller conditioning set $C_{\mathrm{up}}(i,j) \subseteq C(i,j)$. By working with a smaller conditioning set, the effective sample size $n_{C(i,j)} = n - |C(i,j)|$ is increased; in this context see also Wille and Bühlmann (2006). In Sections 5.1 and 5.2, we detail this reasoning for undirected graphs and DAGs, respectively. For bidirected graphs, the conditioning set occurring in the testing problem (3.3) is already as small as possible as it is the empty set $C(i,j) = \varnothing$.

### 5.1 Decreasing the Size of the Conditioning Set in Undirected Graphs

The following graphical condition is the key to finding a smaller conditioning set $C_{\mathrm{up}}(i,j) \subseteq C(i,j)$.

LEMMA 6. *Suppose the observed random vector $Y$ is distributed according to a multivariate normal distribution that is pairwise faithful to an undirected graph $G = (V, E)$. Let $i, j \in V$ be two vertices and define $G^{ij}$ to be the subgraph of $G$ obtained by removing the edge $i$ — $j$, which may or may not be present in $G$. Let $C \subseteq V \setminus \{i, j\}$ be a subset that separates $i$ and $j$ in $G^{ij}$. Then,*

$$Y_i \perp\!\!\!\perp Y_j \mid Y_{V \setminus \{i,j\}} \quad \Longleftrightarrow \quad Y_i \perp\!\!\!\perp Y_j \mid Y_C.$$

PROOF. ($\Longrightarrow$): By the faithfulness assumption, $G$ does not contain the edge $i$ — $j$, so $G = G^{ij}$. Thus the global undirected Markov property for $G$ (see Section 2.1) implies $Y_i \perp\!\!\!\perp Y_j \mid Y_C$.

($\Longleftarrow$): Let $\mathrm{nb}(j)$ be the set of vertices in $V \setminus (\{i, j\} \cup C)$ that, in the graph $G^{ij}$, are connected to $j$ by an edge. In $G$, the set $C \cup \{j\}$ separates $i$ and $\mathrm{nb}(j)$ and thus we obtain via the global Markov property for $G$ that

$$Y_i \perp\!\!\!\perp Y_{\mathrm{nb}(j)} \mid Y_{C \cup \{j\}}. \tag{5.1}$$

Applying standard properties of conditional independence (Lauritzen (1996), Section 3), we obtain from (5.1) and the assumed $Y_i \perp\!\!\!\perp Y_j \mid Y_C$ that

$$Y_i \perp\!\!\!\perp Y_j \mid Y_{C \cup \mathrm{nb}(j)}. \tag{5.2}$$

Moreover, in the graph $G$, the set $C \cup \mathrm{nb}(j) \cup \{i\}$ separates $j$ from the remaining vertices $V \setminus (\{i, j\} \cup C \cup \mathrm{nb}(j))$. Hence, by the global Markov property for $G$,

$$Y_j \perp\!\!\!\perp Y_{V \setminus (\{i,j\} \cup C \cup \mathrm{nb}(j))} \mid Y_{C \cup \mathrm{nb}(j) \cup \{i\}}, \tag{5.3}$$

which in conjunction with (5.2) implies that $Y_i \perp\!\!\!\perp Y_j \mid V \setminus \{i, j\}$. (Note that this proof could be reduced to a single application of the global Markov property if global faithfulness was assumed about the data-generating distribution. Under global faithfulness, $Y_i \perp\!\!\!\perp Y_j \mid Y_C$ implies $G = G^{ij}$.) $\quad\square$

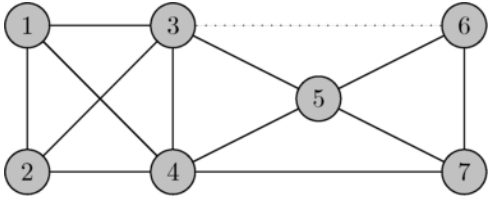

Consider now testing an uncertain edge $i — j$ in $E_u$ by testing the hypothesis $H_{ij} : \rho_{ij \cdot V \setminus \{i,j\}} = 0$. Remove the edge $i — j$ from $G_{up}$ to obtain the graph $G_{up}^{ij}$. In this known graph $G_{up}^{ij}$ we can determine a subset $C_{up}(i, j) \subseteq V \setminus \{i, j\}$ that separates $i$ and $j$. (Choosing this subset to be of minimal cardinality yields the largest gain in effective sample size.) We know a priori that the true data-generating distribution is pairwise faithful to an undirected graph $G$ which is a subgraph of the known graph $G_{up}$. Since graphical separation in $G_{up}$ implies graphical separation in the subgraph $G$, we can deduce from Lemma 6 that

$$(5.4) \qquad \rho_{ij \cdot V \setminus \{i,j\}} = 0 \iff \rho_{ij \cdot C_{up}(i,j)} = 0.$$

As an example, consider as an upper graph the undirected graph in Figure 5. For testing the edge $3 — 6$, which is drawn dotted, we can test $H_{36} : \rho_{36 \cdot 12457} = 0$ but also simply $H_{36} : \rho_{36 \cdot 45} = 0$ or $H_{36} : \rho_{36 \cdot 57} = 0$.

Based on (5.4), the model selection testing problem in (3.3) can be replaced by the problem of testing the $q$ hypotheses

$$(5.5) \qquad H_{ij} : \rho_{ij \cdot C_{up}(i,j)} = 0 \quad \text{vs.} \quad K_{ij} : \rho_{ij \cdot C_{up}(i,j)} \neq 0,$$
$$(i, j) \in E_u.$$

These hypotheses can again be tested using the corresponding sample partial correlations $r_{ij \cdot C_{up}(i,j)}$, or rather the sample $z$-transforms $z_{ij \cdot C_{up}(i,j)}$. Proposition 5 still holds for the vector $(z_{ij \cdot C_{up}(i,j)} \mid (i, j) \in E_u)$. Instead, one could also work with more efficient maximum likelihood estimates computed for the model $N(G_{up})$. It should be noted, however, that the $z$-transform need no longer be variance-stabilizing when applied to such maximum likelihood estimates of partial correlations (Roverato (1996)).

## 5.2 Decreasing the Size of the Conditioning Set in DAGs

In a DAG with well-numbered vertex set, the conditioning set $C(i, j) = \{1, \dots, j\} \setminus \{i, j\}$ may be reduced to any subset $C_{up}(i, j)$ that d-separates $i$ and

$j$ in the graph $G_{up}^{ij}$, defined to be the upper DAG $G_{up}$ with the edge $i \longrightarrow j$ removed. The validity of this replacement can be established using Lemma 7. We note that a simple choice for such a set $C_{up}(i, j)$ are the parents of $j$ in $G_{up}^{ij}$, that is, the set of vertices $k$ that are such that $k \longrightarrow j$ in $G_{up}^{ij}$. However, this need not be the d-separating set in $G_{up}^{ij}$ of smallest cardinality.

Consider, for example, the DAG in Figure 6 as an upper DAG. For testing the edge $1 — 5$, which is drawn dotted, we can use the parents of 5 to test $H_{15} : \rho_{15 \cdot 34} = 0$ but alternatively we can test $H_{15} : \rho_{15 \cdot 2} = 0$.

LEMMA 7. *Suppose the observed random vector $Y$ is distributed according to a multivariate normal distribution that is pairwise faithful to a DAG $G = (V, E)$ with well-numbered vertex set. For any vertices $i, j \in V$ define $G^{ij}$ to be the subgraph of $G$ obtained by removing the edge $i \longrightarrow j$, which may or may not be present in $G$. Let $C \subseteq \{1, \dots, j\} \setminus \{i, j\}$ be a subset that d-separates $i$ and $j$ in $G^{ij}$. Then,*

$$(5.6) \qquad Y_i \perp\!\!\!\perp Y_j \mid Y_{\{1,\dots,j\} \setminus \{i,j\}}$$
$$\iff Y_i \perp\!\!\!\perp Y_j \mid Y_C.$$

PROOF. ($\Longrightarrow$): By the faithfulness assumption, $G$ does not contain the edge $i \longrightarrow j$, so $G = G^{ij}$. Thus the global directed Markov property for $G$ (see Section 2.3) implies $Y_i \perp\!\!\!\perp Y_j \mid Y_C$.

($\Longleftarrow$): Let $\tilde{G} = (\tilde{V}, \tilde{E})$ be the subgraph of $G$ induced by $\{1, \dots, j\}$ and set $\tilde{Y} = Y_{\{1,\dots,j\}}$. Then $\tilde{V}$ is well-numbered for $\tilde{G}$ and $\tilde{Y}$ is pairwise faithful to $\tilde{G}$. Because $C$ d-separates $i$ and $j$ in $\tilde{G}^{ij}$ if and only if $C$ d-separates $i$ and $j$ in $G^{ij}$ and because (5.6) involves only $\tilde{Y}$, we may assume for the proof that $G = \tilde{G}$ and $Y = \tilde{Y}$. Note that $j$ is now a terminal vertex in $G$; that is, $j$ has no children in $G$.

Let pa($j$) be the set of parents of $j$ in $G^{ij}$. We claim that pa($j$) $\setminus C$ and $i$ are d-separated in $G$ given $C$. For, any path $\gamma$ between $i$ and some $k \in$ pa($j$) $\setminus C$ in $G$ either includes $j$ as a nonendpoint or does not.

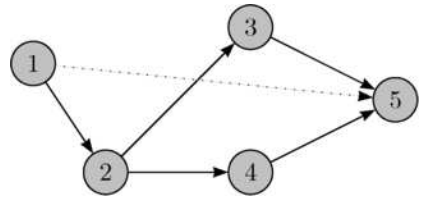

In the first case, since $j$ is terminal in $G$ it must be a collider in $\gamma$ with no directed path to any $\bar{c} \in C$, so $\gamma$ is blocked relative to $C$ in $G$. In the second case, consider the extension of the path $\gamma$ given by appending the edge $k \longrightarrow j$. By the assumed d-separation of $i$ and $j$ in $G^{ij}$ given $C$ and the fact that $k \notin C$ is a noncollider in the extended path, the path $\gamma$ must again be blocked relative to $C$ in $G$. This establishes the claim.

Hence, by the global directed Markov property for $G$,

$$(5.7) \qquad Y_i \perp\!\!\!\perp Y_{\mathrm{pa}(j) \setminus C} \mid Y_C.$$

Because $Y$ has a multivariate normal distribution, it follows from (5.7) and the assumed independence $Y_i \perp\!\!\!\perp Y_j \mid Y_C$ that $Y_i \perp\!\!\!\perp (Y_j, Y_{\mathrm{pa}(j) \setminus C}) \mid Y_C$, hence

$$(5.8) \qquad Y_i \perp\!\!\!\perp Y_j \mid Y_{C \cup \mathrm{pa}(j)}.$$

Since $j$ is terminal in $G$, the set $C \cup \mathrm{pa}(j) \cup \{i\}$ d-separates $j$ from the remaining vertices $\{1, \dots, j\} \setminus (\{i,j\} \cup C \cup \mathrm{pa}(j))$ in $G$. Therefore, by the global Markov property for $G$,

$$(5.9) \qquad Y_j \perp\!\!\!\perp Y_{\{1,\dots,j\} \setminus (\{i,j\} \cup C \cup \mathrm{pa}(j))} \mid Y_{C \cup \mathrm{pa}(j) \cup \{i\}},$$

which in conjunction with (5.8) implies that $Y_i \perp\!\!\!\perp Y_j \mid Y_{\{1,\dots,j\} \setminus \{i,j\}}$. (This proof could be reduced to a single application of the global Markov property for $Y$.) □

REMARK 8. The proof of implication ($\Longleftarrow$) in Lemma 6 holds for any distribution, not necessarily Gaussian, that obeys the global Markov property of the graph $G$. This is in contrast to the proof of implication ($\Longleftarrow$) in Lemma 7, where we have employed a special property of the multivariate normal distribution when deducing (5.8). This special property, namely the fact that $Y_i \perp\!\!\!\perp Y_j$ and $Y_i \perp\!\!\!\perp Y_k$ implies $Y_i \perp\!\!\!\perp (Y_j, Y_k)$, is crucial. For example, it is easy to choose a joint distribution for a binary random vector $(Y_1, Y_2, Y_3)^t$ such that $Y_1 \perp\!\!\!\perp Y_2$ and $Y_1 \perp\!\!\!\perp Y_3$ but $Y_1 \not\perp\!\!\!\perp (Y_2, Y_3)$. This distribution is then pairwise but not globally faithful to the DAG $1 \longrightarrow 3 \longleftarrow 2$, and since $Y_1 \not\perp\!\!\!\perp Y_3 \mid Y_2$, it yields a contradiction to the claim of Lemma 7 for binary random variables. A way around an assumption of global faithfulness is to apply Lemma 6 to an undirected graph $G_{\mathrm{up}}^m$ such that $N(G_{\mathrm{up}}) \subseteq N(G_{\mathrm{up}}^m)$. Such a graph $G_{\mathrm{up}}^m$ can be obtained via the moralization procedure (Lauritzen (1996)).

## 6. DISCUSSION

Gaussian graphical models are determined by pairwise (conditional) independence restrictions, which are in correspondence to the edges that are absent from the underlying graph. These restrictions can be converted into a set of hypotheses that can be tested in order to select a model, or equivalently, a graph. If the arising issue of multiple testing is appropriately addressed, then the selection of a graph can be performed while controlling error rates for incorrect edge inclusion. As reviewed in Section 3, the literature provides a number of methods for such error rate control.

In graphical model selection, controlling incorrect edge inclusion allows us to detect the most important features of multivariate dependence patterns. However, the graph encoding these features need not necessarily yield a model that fits the data well, and if the choice of such a model is the primary focus, then other model selection methods (e.g., the score-based and Bayesian methods discussed in Section 1) may be preferable.

The number of edges in the graph to which the true data-generating distribution is pairwise faithful equals the number of false null hypotheses $H_{ij}$ in (3.3). It would be interesting to adapt existing methods for estimating or bounding this latter number (see, e.g., Meinshausen and Bühlmann (2005)) to make them applicable in graphical model selection. This would allow us to assess the sparseness of the underlying graph by estimating or bounding the number of its edges. Such knowledge could also be used to design more powerful multiple testing-based methods of graph selection in an empirical Bayes framework (see, e.g., Efron et al. (2001)).

In order to associate unique null hypotheses with acyclic directed graphs (DAGs), we restricted ourselves to the situation where a well-numbering of the variables is known a priori. Requiring the knowledge of such a total order is clearly very restrictive. More commonly, time of observation of variables and other considerations provide a priori knowledge in form of a partial order, which allows us to identify ordered blocks of variables. Such blocking strategies appear in many case studies (Caputo, Heinicke and Pigeot (1999), Caputo et al. (2003), Didelez et al. (2002), Mohamed, Diamond and Smith (1998)); compare also Wermuth and Lauritzen (1990). In our illustration of DAG selection in Section 4 the structure of a metabolic pathway yields a partial order among genes that we extended rather arbitrarily to a total order. Hence, a more appropriate analysis might proceed by using the ordered blocks of variables and a generalization of the multiple testing-based model selection we described in order to select a Gaussian chain graph model; see Drton and Perlman (2007) for details on this generalization.

For the class of ancestral graphs of Richardson and Spirtes (2002), multiple testing-based model selection as presented here is less natural. The reason is that in this class there is no distinguished complete graph whose edges could be tested for absence/presence by testing associated conditional independences. If a particular complete ancestral graph is chosen, then a subgraph could be selected similarly as for the other models. The pairwise Markov property (Richardson and Spirtes (2002), page 979) would yield the conditional independences that would have to be tested. However, this procedure could not assure that the selected ancestral graph is maximal (Richardson and Spirtes (2002), page 978). Gaussian models associated with nonmaximal ancestral graphs cannot in general be specified in terms of conditional independence.

Finally, many of the ideas presented here in the framework of Gaussian models carry over to the case of discrete variables or even the mixed case of discrete and continuous variables. Background on the distributional assumptions in the mixed case can be found, for example, in Lauritzen (1996). However, when adapting the methods reviewed here, care must be taken, as conditional independence in the multivariate normal distribution exhibits special properties not shared by other distributional settings. In particular, moving from conditional independence statements between pairs of random variables to ones involving sets of random variables may be valid in a multivariate normal distribution but not in other distributions; compare Remark 8.

## ACKNOWLEDGMENT

We warmly thank Sandrine Dudoit, Thomas Richardson, Juliet Shaffer and Mark van der Laan for helpful comments. This material is based upon work supported by National Science Foundation Grant 0505612.